\documentclass{amsart}
\usepackage{latexsym,epsfig,amssymb,amsmath,amsthm,color,url,bbm}
\usepackage[comma,sort&compress]{natbib}
\usepackage{graphicx,pstricks}
\allowdisplaybreaks 
\numberwithin{equation}{section}

\newtheorem{Lemma}{Lemma}[section]

\newtheorem{Proposition}[Lemma]{Proposition}
\newtheorem{Corollary}[Lemma]{Corollary}

\theoremstyle{definition}
\newtheorem{Remark}{Remark}[section]
\theoremstyle{definition} \newtheorem{Example}{Example}[section]
\newtheorem{Definition}{Definition}[section]

\def\cinP{\stackrel{P}{\to}}

\def\ba{\boldsymbol a}

\def\bx{\boldsymbol x}

\def\by{\boldsymbol y}

\def\K{\mathcal{K}}
\def\E{\mathbb{E}}

\def\R{\mathbb{R}}

\def\F{\mathcal{F}}
\def\Ffiniteorline{{\F_{\text{finite or line}} }}
\def\Finv{F^{\leftarrow}}

\def\snp{\mathcal{S}_n^{'}}
\def\SM{{   {\mathcal{S}_n^{'}}^M}}

\begin{document}

\title[QQ Plots and Random Sets]{QQ Plots, Random Sets\\ and Data from a
  Heavy Tailed Distribution}

\author[B.\ Das]{Bikramjit\ Das}
\address{Bikramjit\ Das\\School of Operations Research and Industrial Engineering\\
Cornell University \\
Ithaca, NY 14853} \email{bd72@cornell.edu}

\author[S.I.\ Resnick]{Sidney I.\ Resnick}
\address{Sidney Resnick\\
School of Operations Research and Industrial Engineering\\
Cornell University \\
Ithaca, NY 14853} \email{sir1@cornell.edu}


\thanks{Sidney Resnick was partially supported by
NSA grant MSPF-05G-049; Bikramjit Das was supported by the NSF training
grant, Graduate and Postdoctoral Training in Probability and Its
Applications, at Cornell University.}

\begin{abstract} {
The QQ plot is a commonly used technique for informally deciding
whether a univariate random sample of size $n$ comes from
a specified distribution $F$. The QQ plot graphs the sample quantiles
against the theoretical quantiles of $F$ and then a visual check is
made to see whether or not the points are close to a straight
line. For a location and scale family of distributions, the intercept
and slope of the straight line provide  estimates for the shift and
scale parameters of the distribution respectively. Here we consider
the set $\mathcal{S}_n$ of points forming the QQ plot as a random closed set in
$\R^2$. We show that under certain regularity conditions on the
distribution $F$,  $\mathcal{S}_n$ converges in probability to a
closed, non-random set.  In the heavy tailed case where
 $1-F$ is a regularly
varying function, a similar result can be shown but a modification is
necessary to provide a statistically sensible result since typically
$F$ is not completely known.}
\end{abstract}
\maketitle

\section{Introduction}\label{sec:intro}
Given a random sample of univariate data
  points, a pertinent question  is whether this sample comes
  from some specified distribution $F$. A variant question is whether
  the sample is from a location/scale family derived from $F$.
Decision techniques are based on  how close the empirical
  distribution of the sample and the distribution $F$ are for some
  sample size $n$.  The empirical distribution function of the iid
  random variables $X_1,\dots,X_n$ is
$$F_n(x):= \frac{1}{n}\sum\limits_{i=1}^n I(X_i \le x), \qquad
  -\infty < x < \infty. $$ The Kolmogorov-Smirnov (KS) statistic is one way
  to measure   the distance between the empirical distribution function
  and the
  distribution function F. Glivenko and Cantelli
  showed (see, for example, \cite{serfling:1980})
  that the KS-statistic converges to 0 almost surely.  The QQ  (or
  quantile-quantile) plot is
  another commonly used device to graphically, quickly
 and informally test the goodness-of-fit of a
  sample in an exploratory way. It has the advantage of being a
  graphical tool, which is
  visually appealing and  easy to understand.  The QQ plot
measures how close the sample quantiles are to the theoretical quantiles.
 For $0 < p< 1 $, the $p^{th}$ \emph{ quantile} of $F $ is
  defined by
  \begin{align} \label{eqn:defFinv}
     F^{\leftarrow}(p) & := \inf\{x: F(x) \ge p\}.
  \end{align}
  The sample $p^{th}$ quantile can be similarly defined as
  $F_n^{\leftarrow}(p)$. If $ X_{1:n} \le
X_{2:n} \le \ldots \le X_{n:n}$ are the order statistics from the
sample, then $F_n^{\leftarrow}(p) = X_{\lceil np
\rceil:n}$, where as usual $\lceil np
\rceil$ is the smallest integer greater than or equal to $np$.
For $0 < p < 1$, $X_{\lceil np \rceil:n}$
is a strongly consistent estimator of $F^{\leftarrow}(p)$
\citep[page 75]{serfling:1980}.

Rather than considering individual quantiles, the QQ plot considers
the sample as a whole and plots the sample quantiles against the
theoretical quantiles of the specified target distribution $F$.  If we
have a correct target distribution, the QQ plot hugs a straight line
through the the origin at an angle of $45^{\circ}$.
Sometimes we have a location and scale family correctly specified up
to unspecified location and scale and in such cases, the
QQ plot concentrates around a straight line with some slope (not
necessarily $45^{\circ}$) and intercept (not necessarily $0$); the
slope and intercept estimate the scale and location.
 Since a variety of
estimation and inferential procedures in the practice of
statistics depends on the
assumption of normality of the data, the normal QQ plot is one of the
most commonly used.

Our goal here is to formally prove the convergence of the QQ plot
(perhaps suitably modified) to a straight line. This would show the
asymptotic consistency of the QQ plot.  The QQ plot formed by a
sample of size $n$ can  be considered  a closed subset of $\R^2$
denoted by $\mathcal{S}_n$. This  set of points that form the QQ
plot in $\R^2$ is
   \begin{align}\label{eqn:defSn}
     \mathcal{S}_n & :=\{(F^{\leftarrow}(\frac{i}{n+1}),X_{i:n}),~~ 1
     \le i \le n
     \}
     \end{align}
where the function $F^{\leftarrow}(\cdot)$ is defined by
 (\ref{eqn:defFinv}).
       For each $n$, $\mathcal{S}_n $ is a random closed set. Note that, if $\{ \mathcal{S}_n\}$  has an
    almost sure limit $S$ then this limit set by the
Hewitt-Savage $0-1$ law
    must be almost surely  constant.
     A straight line (or some closed subset of a straight line) is also a
    closed set in $\R^2$.  Under certain regularity conditions on $F$,
    we show that the random set $\mathcal{S}_n$
    converges in probability to a straight line (or some closed subset
     of a straight line), in a suitable topology on closed subsets of
 $\R^2$.

Section \ref{sec:prelim} is devoted to preliminary results on the
convergence of random closed sets.  We also discuss a result on
convergence of quantiles and, because of our interest in heavy tails,
we introduce the concept of regular variation. In Section
\ref{sec:qqplot}, we assume the random variables have a specified
distribution $F$ and we consider convergence of the random closed sets
$ \mathcal{S}_n$ forming the QQ plot. In Section \ref{sec:qqplotrv},
the idea of the QQ plot is extended to the case where we know that the
data is heavy tailed, that is $1-F$ is regularly varying. We assume we
do not know the exact distribution of $F$; we presume the distribution
is heavy tailed but do not know either the tail index or the slowly
varying component. The usual QQ plot is not informative in a
statistical sense and hence must be modified by a thresholding
technique.

In Corollary \ref{cor:par} we
   have convergence of a log-transformed version of the QQ plot
   to a straight line when the distribution of the random sample is
   Pareto. Now Pareto being a special case of a distribution with
   regularly varying tail, we use the same plotting technique for
   random variables having a regularly varying tail after thresholding
   the data.  We provide
   a convergence in probability result considering the $k=k(n)$ upper order
   statistics of the data set where $k \rightarrow \infty$ and $k/n
   \rightarrow 0$.  In Section \ref{sec:lsline}, a continuity result
   is provided for a least squares line through these special kinds of
   closed sets. See \citet{kratz:resnick:1996,
beirlant:vynckier:teugels:1996}.

\section{Preliminaries}\label{sec:prelim}

\subsection{Closed sets and the Fell topology} \label{subsec:convran}
We denote the distance between the points $\bx$ and $\by$ by
$d(\bx,\by)$; $\mathcal{F}, \mathcal{G}$ and $ \mathcal{K} $ are the
classes of closed, open and compact subsets of $\mathbb{R}^d$
respectively. These quantities are sometimes subscripted by the
dimension of the space if this needs to be emphasized for clarity. We
are interested in closed sets because the sets of
interest such as $\mathcal{S}_n$ are random closed sets.  There are
several ways to define a topology on the space of closed sets. The
Vietoris topology and the Fell topology are frequently used and these
are hit-or-miss kinds of topologies. We shall discuss the Fell
topology below. For further discussion refer to \citet{beer:1993,
matheron:1975, molchanov:2005}.

For a set $ B \subset \mathbb{R}^d $,  define $\mathcal{F}_B$ as the
class of closed sets
hitting $B$ and $\mathcal{F}^B$ as the class of closed sets disjoint
from $B$:
$$
   \mathcal{F}_B = \{ F: F \in \mathcal{F}, ~~ F \bigcap B \neq \emptyset
  \}, \quad
 \mathcal{F}^B = \{ F: F \in \mathcal{F}, ~~ F \bigcap B =
\emptyset
  \}.
$$
Now the space $\mathcal{F}$ can be topologized by
the Fell topology which has as its subbase the families
$\{\mathcal{F}^K, K \in \mathcal{K}\}$ and $\{\mathcal{F}_G, G \in
\mathcal{G}\}$.

  A sequence $\{F_n\}$ converges in the Fell topology
 towards a limit $F$ in $\F$ (written
$F_n \to F$) if and only if it satisfies two conditions:
\begin{enumerate}
\item If an open set $G$ hits $F$, $G$ hits all  $F_n$, provided $n$
  is sufficiently large.
 \item If a compact set $K$ is
disjoint from $F$, it is disjoint from $F_n$ for all sufficiently
large $n$.
\end{enumerate}
The following result \citep{matheron:1975}
 provides useful conditions for convergence.

\begin{Lemma}\label{lemma:conv1}
     For $F_n,F \in \mathcal{F}, n \ge 1,\, F_n \rightarrow F$ as
     $n\rightarrow \infty$ if and only if the following two conditions
     hold

\begin{align}
&\ \text{For any $\by \in F$,
for all large $n$,
there exists  $\by_n \in F_n$ such that $d(\by_n,\by) \rightarrow 0$ as
      $n \rightarrow \infty$.}\label{eqn:crit1}\\
&\ \text{For any subsequence $\{ n_k\}$, if $\by_{n_k} \in F_{n_k}$
      converges, then $\lim_{k \rightarrow \infty} \by_{n_k} \in
      F$.}\label{eqn:crit2}
\end{align}
      Furthermore, convergence of sets $\mathcal{S}_n \rightarrow
       \mathcal{S}$ in $\mathcal{K}$ is equivalent to the analogues
       of (\ref{eqn:crit1}) and (\ref{eqn:crit2}) holding as well as
       $\sup_{j \ge 1} \sup\{\|\bx \|:\bx \in S_j \} < \infty$ for some
       norm $\|\cdot \|$ on $\mathbb{R}^d$.

     \end{Lemma}

      Note that if the sets are random elements of $\mathcal{K}$
      and $S \in \mathcal{K}$ is non-random, then Lemma \ref{lemma:conv1}
      can be used to characterize almost sure convergence or
      convergence in probability. We are going to define random sets
      in the next subsection.
\begin{Definition}[Hausdorff Metric]
        Suppose $d:\R^d\times \R^d \to \R_+$ is a metric on
  $\mathbb{R}^d$. Then for $S,T \in \mathcal{K}$, define the
  Hausdorff metric \citep{matheron:1975} $D:\mathcal{K} \times
  \mathcal{K} \rightarrow \mathbb{R}_{+}
  $ by
\begin{align}
D(S,T) & = \inf\{ \delta: S \subset T^{\delta},T \subset
     S^{\delta}\} , \label{eqn:defD}\\
  \intertext{ where for $S \in \mathcal{K}$ and $\delta > 0$, }
    S^{\delta} & = \{\bx: d(\bx,\by)< \delta \mathrm{~~for~~ some~~} \by \in
     S\} \label{eqn:defS}
   \end{align}
is the $\delta$-neighborhood or $\delta$-swelling of $S$.
   \end{Definition}
The topology usually used on  $\mathcal{K}$ is the {\it myopic
  topology\/} with sub-base elements $\{\mathcal{K}^F, F \in
  \mathcal{F}\}$ and  $\{\mathcal{K}_G, G \in
  \mathcal{G}\}.$  The myopic topology on $\mathcal{K}$ is stronger than
  the Fell topology relativized to $\mathcal{K}$. The topology on
  $\mathcal{K}'= \mathcal{K}\setminus \{\emptyset\}$ generated by the
  Hausdorff metric  is equivalent to the myopic topology on $\mathcal{K}'$
\citep[page 405]{molchanov:2005}.

In certain cases, convergence on $\F$ can be reduced to convergence on
$\K$.
\begin{Lemma}\label{lem:reduceConvToK}
Suppose $F_n$ and $F$ are closed sets in $\F$ and
and that there exist $\K_1 \subset \K$ satisfying
\begin{enumerate}
\item $\displaystyle \bigcup_{K\in \K_1} K =\mathbb{E}.$
\item For $\delta >0$ and $K \in \K$, we have
$\overline{K^\delta} \in \K_1.$
\item $\displaystyle F_n \bigcap K \to F\bigcap K,\quad \forall K \in
  \K_1.$
\end{enumerate}
Then $F_n \to F$ in $\F$.
\end{Lemma}

\begin{Remark}\label{rem:convFalse}
The converse is false. Let $\mathbb{E}=\mathbb{R}$, $F_n=\{1/n\}$,
$F=\{0\}$ and $K=[-1,0]$. Then $F_n \to F$ but
$$F_n \bigcap K =\emptyset \not\to F \bigcap K=F.$$
The operation of  intersection is not a continuous operation in
$\F\times \F$
\citep[page 400]{molchanov:2005}; it is only upper semicontinuous
\citep[page 9]{matheron:1975}.
\end{Remark}

\begin{proof} We use Lemma \ref{lemma:conv1}. If $x \in F,$
there exists $K \in \K_1$ and $x\in K$. So $x\in F\cap K$ and from Lemma
\ref{lemma:conv1}, since $F_n\cap K \to F\cap K$ as $n \to
\infty$,
 we have existence of
$x_n \in F_n \cap K$ and $x_n\to x$. So we have produced $x_n \in F$
and $x_n \to x$ as required for \eqref{eqn:crit1}.

To verify \eqref{eqn:crit2},  suppose $\{x_{n_k} \}$ is a subsequence
such that
$x_{n_k} \in F_{n_k}$ and $\{x_{n_k} \}$ converges to, say,
$x_\infty$. We need to show $x_\infty \in F$. There exists $K_\infty
\in \K_1$ such that $x_\infty \in K_\infty$. For any $\delta >0$,
$x_{n_k} \in \overline{ K^\delta_\infty } \in \K_1$ for all sufficiently large
$n_k$. So $x_{n_k}\in F_{n_k}\cap \overline{ K^\delta_\infty } .$ Since
$F_{n_k}\cap \overline{ K^\delta_\infty } \to F\cap \overline{K^\delta_\infty },$
we have
$\lim_{k \to \infty} x_{n_k}=x_\infty \in F\cap \overline{K^\delta_\infty }.$
So $x \in F$.
\end{proof}

The next result shows when a point set approximating a curve actually
converges to the curve. For this Lemma, $C(0,1]$ is the class of real
  valued continuous functions on $(0,1]$ and $D_l(0,\infty]$ is the
      class of left continuous functions on $(0,\infty]$ with finite
    right hand limits.

\begin{Lemma}\label{lem:approxCurve}
Suppose $0\leq x(\cdot) \in C(0,1]$ is continuous on $(0,1]$ and strictly
    decreasing with $\lim_{\epsilon \downarrow 0}
    x(\epsilon)=\infty.$  Suppose further that $y_n(\cdot) \in
    D_l(0,1]$ and $y(\cdot) \in C(0,1]$ and $y_n \to y$ locally
    uniformly on $(0,1]$; that is, uniformly on compact subintervals
    bounded away from $0$. Then for $k=k(n) \to \infty$,
$$F_n:=\{\bigl(x(\frac {j}{k} ),y_n(\frac {j}{k} )\bigr) ; 1 \leq j \leq k\}
\to F:=\{\bigl(x(t),y(t)\bigr) ; 0<t\leq 1\}
=\{\bigl(u,y(x^\leftarrow (u))\bigr); x(1)\leq u <\infty\},$$
in $\F$.
\end{Lemma}

\begin{proof}
Pick $t\in (0,1]$, so that $(x(t),y(t))\in F$. Then
$$F_n \ni \bigl( x( \lceil kt \rceil /k),y_n( \lceil kt \rceil /k)\bigr)
\to \bigl( x(t),y(t)\bigr)\in F,$$
in $\mathbb{R}^2,$  verifying \eqref{eqn:crit1}. For
\eqref{eqn:crit2}, Suppose
$\bigl( x(j(n')/k(n'), y_{n'}(j(n')/k(n')\bigr) \in F_{n'}$ is a
convergent subsequence in $\mathbb{R}^2$. Then $\{ x(j(n')/k(n')\}$ is
convergent in $\mathbb{R}$ and because $x(\cdot)$ is strictly
monotone,
$\{j(n')/k(n')\} $ converges to some $l\in (0,1]$. Then
$$ F_{n'} \ni  \bigl( x(j(n')/k(n'), y_{n'}(j(n')/k(n')\bigr)
\to \bigl(x(l),y(l)
\bigr) \in F,$$
which verifies \eqref{eqn:crit2}.
\end{proof}

\subsection{Random closed sets and weak convergence}
In this section, we review definitions and characterizations of weak
     convergence of random closed sets.  In subsequent sections we
     will show convergence in probability, but since the limit sets
     will be non-random, weak convergence and convergence in
     probability coincide.  See also \citet{matheron:1975,
     molchanov:2005}.

Let $( \Omega, \mathcal{A}, P^{'})$ be a complete probability
space. $\F$ is the space of all closed sets in $\R^d$
    topologized by the Fell topology.
    Let  $\sigma_{\F}$ denote the Borel $\sigma $-algebra generated by
    the Fell topology of open sets.
A \emph{random closed set} $X: \Omega \mapsto \F $ is  a measurable
mapping from $(\Omega,\mathcal{A},P^{'})$ to $ (\F, \sigma_{\F}
    )$.  Denote by $P$ the induced probability on $\sigma_{\F}$,
    that is, $P=P^{'}\circ X^{-1}.$
A sequence of  random closed sets $\{X_n\}_{n\ge 1} $ weakly
  converges to a  random closed set $X$ with distribution $P$) if the
  corresponding induced
  probability measures $\{P_n\}_{n\ge 1} $ converge weakly to $P$,
  i.e.,
  \[ P_n(\mathcal{B})=P_n^{'}\circ X_n^{-1}( \mathcal{B})
\to P(\mathcal{B})=
P^{'}\circ X^{-1}( \mathcal{B}),
 \qquad \text{ as } n \to \infty, \]
for each $ \mathcal{B} \in \sigma_{\mathcal{F}} $ such that
$P(\partial \mathcal{B}) =
  0$.

This is not always straightforward to verify from the definition.
We find useful the following characterization of weak convergence in terms of
sup-measures \citep{vervaat:1997}.
Suppose $h: \R^d \mapsto \R_{+} =[0,\infty)$.
For $X\subset \mathbb{R}^d$, define
$h(X) = \{ h(x): x \in X \}$ and $h^\vee $ is the sup-measure
generated  by $h$
defined by
$$h^\vee (X)=\sup\{ h(x): x \in X \}$$
 \citep{molchanov:2005, vervaat:1997}.
These definitions
permit the following  characterization
 \citep[page 87]{molchanov:2005}.
     \begin{Lemma}  \label{thm:choq}
     A sequence $\{X_n\}_{n\ge 1}$ of random closed sets converges weakly to
     a random closed set $X$ if and only if $\E h^{\vee} (X_n)$
     converges to $\E h^{\vee} (X)$ for every non-negative continuous
     function $h:\R^d \mapsto \R$ with a bounded support.
     \end{Lemma}

\subsection{Convergence of sample quantiles}\label{subsec:conqn}
The sample quantile is a strongly consistent estimator
of the population quantile (\cite{serfling:1980}, page 75). The weak
consistency of sample quantiles as estimators of
 population quantiles was shown by
\cite{smirnov:1949}; see also \cite[page 179]{resnick:1999book}.
We will make use of the Glivenko-Cantelli lemma describing uniform
convergence of the sample empirical distribution and also take note of
the following quantile estimation result.
\begin{Lemma} \label{lemma:quantile}
Suppose $F$ is strictly increasing at $F^{\leftarrow}(p)$ which
means that for all $\epsilon > 0$,
\begin{equation*}
 F(F^{\leftarrow}(p - \epsilon)) < p < F(F^{\leftarrow} (p + \epsilon)).
\end{equation*}
Then we have that the $p^{th}$ sample quantile, $X_{\lceil
np\rceil:n}$ is a weakly consistent quantile estimator,
$$ X_{\lceil np\rceil:n} \cinP F^{\leftarrow}(p)$$
As before, $ \lceil np \rceil$ is the  $1^{st} $ integer $\ge
np$ and $X_{i:n}$ is the $i^{th}$ smallest order statistic. \end{Lemma}

\subsection{Regular variation} Regular variation is the mathematical
underpinning of  heavy tail analysis. It is discussed in many books
such as \cite{resnickbook:2006, resnick:1987, seneta:1976,
geluk:dehaan:1987, dehaan:1970, dehaan:ferreira:2006,
bingham:goldie:teugels:1987}.

\begin{Definition}[Regular variation]\label{subsec:rv}
     A measurable function $U(\cdot): \R_+ \to \R_+$ is  regularly
      varying at $\infty$ with  index $\rho \in \R$ if for $x > 0$
      \begin{equation} \label{eqn:defrv1}
      \lim\limits_{t\to \infty} \frac{U(tx)}{U(t)} = x^{\rho}.
      \end{equation}
      We write $U \in RV_{\rho}$.
     \end{Definition}
\begin{Remark}
When $\rho=0 $ we call $U(\cdot)$ \textit{slowly varying } and
       denote it by $L(\cdot)$. For $\rho \in \R $,  we can always
       write $U  \in RV_{\rho}$ as:
      \begin{equation} \label{eqn:defrv2}
      U(x) = x^{\rho} L(x)
      \end{equation}
where $L(\cdot)$ is slowly varying.
\end{Remark}

\section{QQ plots from a known distribution: Random sets  converging
  to a constant set}
      \label{sec:qqplot}
  In this section, we will use the results in Section
\ref{sec:prelim} to  show the convergence of the random closed
sets given by \eqref{eqn:defSn} consisting of the points forming the
QQ plot to a non-random  set in $\R^2$.  First we consider the easiest
case where the random variables are iid from a uniform distribution.
Then we consider more general distributions which are continuous and
strictly increasing on their support. This result will be derived from
the uniform case. Because we are interested in heavy tailed
distributions, our final corollary in this section is about the Pareto
distribution which is the exemplar of the heavy tailed distribution.

\subsection{The Uniform case}\label{subsec:qqunif}
The first simple example is QQ plot from the uniform distribution.
\begin{Proposition}\label{propn:unif}
Suppose $U_1,U_2,\ldots,U_n$ are iid U(0,1). Denote the order
       statistics of this sample by $ U_{1:n} \le U_{2:n}
       \le \ldots \le U_{n:n}$. Define
       \begin{align}
     \mathcal{S}_n & :=\{(\frac{i}{n+1},U_{i:n}),~~ 1 \le i \le n \}
     \label{eqn:defSnU}\\
  \intertext{and}
     \mathcal{S} & :=\{(x,x): 0 \le x \le 1\} \label{eqn:defSU}.
   \end{align}
   Then $\mathcal{S}_n \stackrel{a.s.}{\to} \mathcal{S}$ in $\K_2$.
       \end{Proposition}

\begin{proof} We apply the convergence criterion
 given in Lemma \ref{lemma:conv1}.
The empirical distribution $U_n(x)=n^{-1}\sum_{i=1}^n I(U_i\leq x) $
converges uniformly for almost all sample paths to $x,\,0\leq x\leq
1$. Without loss of generality suppose this true for all sample
paths. Then for all sample paths,
 the same is true  for the inverse process $U_n^\leftarrow
(p)=U_{\lceil np\rceil :n},\,0\leq p\leq 1$; that is
$$\sup_{0\leq p\leq 1} |U_{\lceil np\rceil :n} -p|\to 0, \quad (n\to\infty).$$

Pick $0\leq y\leq 1$ and let $ \by=(y,y) \in \mathcal{S}$.  For each n,
define $\by_n$ by
  \begin{equation}\label{eqn:defyn}
            \by_n = \Bigl(\frac{\lceil ny \rceil}{n+1},U_{\lceil ny\rceil
            :n}\Bigr),
           \end{equation}
so that $\by_n \in \mathcal{S}_n$.  Since $ |ny - \lceil ny
          \rceil| \le 1$ ,
 $ {\lceil ny \rceil}/{n+1} \to y$ and since $U_{\lceil ny\rceil
            :n}\to y$, we have $\by_n \to (y,y) \in \mathcal{S}$.
Hence criterion (\ref{eqn:crit1}) from Lemma
 \ref{lemma:conv1} is satisfied.

Now suppose we have a subsequence $\{n_k\}$ such that
$\by_{n_k} \in \mathcal{S}_{n_k}$
converges.  Then $\by_{n_k}$ is of the form
$\by_{n_k}=(i_{n_k}/(n_k+1),U_{i_{n_k}:n_k})$ for some $1 \le i_{n_k} \le
n$ and  for some $x \in [0,1],$ we have
 $i_{n_k}/(n_k+1) \to x$ and hence also
 $i_{n_k}/n_k \to x$. This implies
$$U_{i_{n_k}:n_k}=U_{\lceil  n_k \cdot \frac{i_{n_k}}{n_k} \rceil :
  n_k} \to x,$$
and therefore $\by_{n_k} \to (x,x)$ as required for
\eqref{eqn:crit2}.\end{proof}

\subsection{Convergence for more general distributions} \label{subsec:qqgen}
Now consider a distribution function $F$ which is more general than
 the uniform, assuming  that $F$ is strictly increasing and continuous on its
 support so
 that $F^{\leftarrow}$ is unique.

\begin{Proposition}\label{propn:genF}
 Suppose $X_1,\ldots, X_n$ are iid with common distribution
$F(\cdot)$ and $ X_{1:n} \le X_{2:n} \le \ldots \le X_{n:n}$
 are the order statistics from this sample.
   If $F$  is strictly increasing and continuous on its support, then
\begin{align*}
 \mathcal{T}_n &: =\{(F^{\leftarrow}(\frac{i}{n+1}), X_{i:n}) ; 1 \le i
\le n\}\\
\intertext{ converges in probability to}
\mathcal{T} &:= \{ (x,x) ; x \in support(F)
\}
\end{align*}
 in $\F_2$.
 \end{Proposition}

 \begin{proof}
According to Lemma \ref{thm:choq}, we must prove for any
non-negative continuous $h:\mathbb{R}^2\mapsto \mathbb{R}_+$ with
compact support
that as $n \to \infty$,
$$\E \bigl( h^\vee (\mathcal{T}_n)\bigr) \to \E \bigl(
h^\vee (\mathcal{T})\bigr).$$

Since $F$ is  continuous,
   $F(X_1), F(X_2), \ldots,F(X_n)$ are iid and uniformly distributed
  on $[0,1]$.   Therefore from Proposition \ref{propn:unif} we have
  that
 \begin{equation} \label{eqn:SntoS}
 \mathcal{S}_n :=\{(\frac{i}{n+1}, F(X_{i:n})) ; 1 \le i \le n\}
\stackrel{d}{=}\{\frac{i}{n+1}, U_{i:n}) ; 1 \le i \le n\}
\stackrel{a.s.}{\to}
\mathcal{S} =\{ (x,x) ; 0 \le x \le 1 \}
 \end{equation}
 in $\K_2$.

We now proceed by considering cases which depend on the nature of the
support of $F$.  We will need the following identity.
For any closed set $X$, function $f:\mathbb{R}^2\mapsto \mathbb{R}_+$
 and function $ \psi : \R^2 \mapsto \R^2 $, we
 have,
 \begin{equation}
  f^{\vee} \circ \psi  (X)  = \sup_{t  \in \psi (X)} f(t)
                       = \sup_{s \in X} f(\psi (s))
                        = \sup_{s \in X} f \circ \psi  (s)
 = (f \circ \psi )^{\vee} (X).\label{eqn:cute}
 \end{equation}

\subsubsection*{Case 1: The support of $F$ is compact, say
$[a,b]$.} This implies $\Finv(0)=a, \,\Finv(1) =b. $
Define the map
 $g:[0,1]^2 \mapsto
[a,b]^2$ by
\[ g(x,y) = (\Finv(x),\Finv(y)). \]
Since $F$ is strictly increasing,
observe that $g(\mathcal{S}_n){=} \mathcal{T}_n$ and
$g(\mathcal{S})= \mathcal{T} $.
Define $ g^{*}: \R^2
\mapsto \R^2$ as the extension of $g$ to all of $\mathbb{R}^2$:
\begin{align*}
g^{*}(x,y) & = (g_1(x), g_1(y)) \\
\intertext{ where }
g_1(z) & = \begin{cases} \Finv(z), & 0 \le z \le 1\\
                         a, & z \le 0\\
                         b, & z\ge 1.
            \end{cases}
 \end{align*}
This makes  $ g^{*}: \R^2 \mapsto \R^2$   continuous.  Since
both $\mathcal{S}_n $ and $\mathcal{S} $ are subsets of $ [0,1] \times
 [0,1]$, we have $g(\mathcal{S}_n) = g^{*}(\mathcal{S}_n)$ and
 $g(\mathcal{S}) = g^{*}(\mathcal{S})$.
 Let $f$ be a continuous function on $\mathbb{R}^2$  with bounded
 support and we have, as
 $n \to \infty$, using \eqref{eqn:cute},
  \begin{align*}
  \E  f^{\vee}  (\mathcal{T}_n) & = \E  f^{\vee} ( g (\mathcal{S}_n))
 = \E f^{\vee} ( g^{*}(\mathcal{S}_n))\\
                      & = \E (f \circ g^{*})^{\vee} (\mathcal{S}_n)
\to \E (f \circ g^{*})^{\vee} (S). \\
\intertext{The previous convergence results from
 $f\circ g^{*}: \R^2 \mapsto \R_+$
being  continuous with bounded support,  $\mathcal{S}_n
 \stackrel{P}{\to} \mathcal{S}$, and  Lemma
 \ref{thm:choq}. The term to the right of the convergence arrow above equals }
                      & = \E f^{\vee} (g^{*}(\mathcal{S}))
                       = \E f^{\vee} (g(\mathcal{S})) = \E f^{\vee}
 (\mathcal{T}).
  \end{align*}
 Therefore   $\mathcal{T}_n$ converges to $T$ weakly and since $T$ is
 a  non-random set, this convergence is  also true in
 probability.

\subsubsection*{Case 2: The support of $F$ is $\R= (-\infty,\infty)$.}
Now define $g:(0,1)^2 \mapsto
\R^2$ by
\[ g(x,y) = (\Finv(x),\Finv(y)). \]
Since $F$ is strictly increasing, $g(\mathcal{S}_n)= \mathcal{T}_n$
 and $g(\mathcal{S}\cap (0,1)^2)= \mathcal{T}$.  Let $f$ be a
 continuous function with bounded support in $[-M,M]^2$, for some
 $M>0$. Extend the definition of $g$ to all of $\R^2$ by
defining $ g^{*}: \R^2 \mapsto \R^2$ as
$$
g^{*}(x,y)  = (g_1(x), g_1(y)),
 $$
where
$$
g_1(z)  = \begin{cases}
\Finv(z), & -M \le \Finv (z) \le M,\\
\Finv(-M),&  \Finv(z) \leq -M,\\
\Finv(M), &  \Finv(z) \geq M.
            \end{cases}
$$
 Therefore $ g^{*}: \R^2 \mapsto \R^2$  is continuous. Now note that
 since $support(f) \subseteq [-M,M]^2$ and $g(x,y) = g^{*}(x,y)$ for
 $(x,y) \in [-M,M]^2$, we will have $f\circ g = f\circ g^{*}$. Therefore
  \begin{align*}
  \E  f^{\vee}  (\mathcal{T}_n) & = \E  f^{\vee} ( g (\mathcal{S}_n))
 = \E (f \circ  g)^{\vee}(\mathcal{S}_n)
                      = \E (f \circ  g^{*})^{\vee}(\mathcal{S}_n)
  \to \E (f \circ g^{*})^{\vee} (S).\\
\intertext{As with Case 1, the convergence follows from
 $f\circ g^{*}: \R^2 \mapsto \R_+$ being continuous with bounded
 support,
$\mathcal{S}_n \stackrel{P}{\to}\mathcal{S}$ and  Choquet's theorem
 \ref{thm:choq}. The term to the right of the convergence arrow equals}
                      & = \E (f \circ g)^{\vee} (S)
 = \E f^{\vee} (g(S)) = \E f^{\vee} (T).
  \end{align*}
 Therefore   $\mathcal{T}_n$ converges to $\mathcal{T}$ weakly. But
 since
$\mathcal{T}$ is a
 non-random set, this convergence is true also in
 probability.

\subsubsection*{Case 3: The support of $F$ is of the form
$[a,\infty)$ or $(-\infty,b]$.}
 This case can be examined  in a similar manner as we have done for Cases 1
 and 2 by considering each end-point of the interval of support of
 $F$ according to its nature.
\end{proof}

\begin{Corollary} \label{cor:exp}
If $F$ is exponential with parameter $\alpha > 0$, i.e., $F(x) = 1-
e^{-\alpha x},\, x> 0$, we have
$$ \{(-\frac{1}{\alpha}\log(1-\frac{i}{n+1}), X_{i:n}); 1 \le i \le
n\} \cinP \{(x,x): 0 \le x < \infty \}.$$
\end{Corollary}

\begin{Corollary} \label{cor:par}
If $F$ is Pareto with parameter $\alpha > 0$, i.e., $F(x) = 1-
x^{-\alpha}$, $x>1$, we have
$$ \{(-\log(1-\frac{i}{n+1}), \log X_{i:n}); 1 \le i \le n\} \cinP
\{( x,\frac{x}{\alpha}): 0 \le x < \infty \}.$$
\end{Corollary}

\section{QQ plots: Convergence of random sets in the regularly varying case }
      \label{sec:qqplotrv}
The classical QQ plot can be graphed only if we know the
target distribution $F$ at least up to location and scale.
We would like to extend the idea of QQ
  plots to the case where the data is from a heavy tailed
  distribution; this is  a semi-parametric assumption
which is more general than assuming the target distribution $F$ is
known up to location and scale.

We model a
  one-dimensional heavy-tailed
  distribution function $F$ by assuming it has a regularly varying
  tail with some   index $-\alpha$, for $\alpha > 0$; that is, if $X$
  has distribution $F$ then,
  \begin{equation}\label{eqn:defRV}
P[X>x] = 1- F(x) = \bar{F}(x) = x^{-\alpha} L(x) , \qquad x>0
\end{equation}
  where L is slowly varying.   In at least an exploratory context,
how   can the QQ plot be used  to validate
this assumption and also to
  estimate $\alpha$? (See \citet[page 106]{resnickbook:2006}.)

 Notice that if we take $L\equiv 1$, $F$ turns out to be a Pareto
  distribution with parameter $\alpha$.  In Corollary (\ref{cor:par}),
  we have seen that if $F$ has a Pareto distribution with parameter
  $\alpha$, then $\mathcal{S}_n$ defined as:
\begin{align}
 \mathcal{S}_n &:= \{(-\log(1-\frac{i}{n+1}), \log X_{i:n}); \; 1 \le
 i \le n\} \label{eqn:Snpar}\\
\intertext{ converges in probability to the set }
   \mathcal{S} &= \{( x, \frac{x}{\alpha})  ;  0 \le x < \infty \}.
\end{align}

 Keeping this in mind, when we have a general $\bar F \in RV_{-\alpha}$,
 let us define $\mathcal{S}_n$ exactly as in (\ref{eqn:Snpar}). Then we are
 able to show that, $\mathcal{S}_n$ converges in probability to the set
\begin{equation} \label{eqn:defSL}
\mathcal{S} = \{(\alpha x, x + \frac{1}{\alpha}\log L(F^{\leftarrow
       }(1-e^{-\alpha x}))) ; \; 0 \le x < \infty
       \}.
\end{equation}
But, since we do not know the slowly varying function
$L(\cdot), $ this result is not  useful for
inference purposes. Estimating $\alpha$ from such a set is not
possible unless $L(\cdot)$ is known,  nor is it clear how
$\mathcal{S}_n $ graphically
approximating such a set would allow us to validate the model
assumption of a regularly varying tail.

Consequently we concentrate on a different asymptotic regime where the
 asymptotic behavior of the random closed set can be freed from
 $L(\cdot)$. For a sample of size $n$ from the distribution $F$, where
 $\bar F \in RV_{-\alpha}$, we consider the upper $k=k(n)$ order statistics of
 the sample where $k(n)/n \to 0$ and construct a QQ plot
 similar to (\ref{eqn:Snpar}).  We assume that $d_\F (\cdot,\cdot)$ is
 some metric on $\F$ which is compatible with the Fell
 topology. Note \cite{flachsmeyer:1964}
 characterized the metrizability of the Fell
 topology and since $\R^d$ is locally compact, Hausdorff and second
 countable his results apply and allow the conclusion that
$\F$ is metrizable under the Fell topology.

For what follows, when $A \in \F_2$, we write
$A+(t_1,t_2)=\{\boldsymbol{a}+(t_1,t_2) : \boldsymbol{a}\in A\}$ for
the translation of $A$.

\begin{Proposition}\label{propn:conrv}
  Suppose we have a random sample
$ X_1,X_2,\ldots,X_n$  from $F$ where $ \bar F \in
RV_{-\alpha}$ and $ X_{(1)} \ge X_{(2)} \ge \ldots \ge X_{(n)}$ are
the order statistics in decreasing order.
 Define
 \begin{align*}
 \mathcal{S}_n & = \{ (-\log{\frac{j}{n+1}},\log X_{(j)}) ; j =1,\ldots,k\}
 \\
  \intertext{ where $k=k(n) \to \infty$ and $k/n \to 0$ as $n
 \to \infty$.
 Also define }
  \mathcal{T}_n & = \{(x,\frac{x}{\alpha}) ; x \ge 0 \} +
(-\log \frac{k}{n+1},\log X_{(k)})
  \\
\intertext{ Then as $n \to \infty$ }
 & \ d_\F (\mathcal{S}_n, \mathcal{T}_n)  \stackrel{P} {\to} 0
 \end{align*}
\end{Proposition}

\begin{Remark}
So after a logarithmic transformation of the data, we make the QQ
plot by only comparing the $k$ largest order
statistics with the corresponding theoretical exponential distribution
quantiles. This  produces an asymptotically linear
plot of slope $1/\alpha$ starting from the point
$(-\log \frac{k}{n+1},\log X_{(k)})$.
\end{Remark}

\begin{proof}
Define
$$
  \mathcal{S}_n^{'}  =\{(-\log \frac{j}{k},
  \log\frac{X_{(j)}}{X_{(k)}}); 1 \le j \le
k\}, \text{ and }
  \mathcal{T}  =\{ (x, \frac{x}{\alpha}); 0 \le x < \infty\}.
$$
 Note that we can write
$$
\mathcal{S}_n^{'}
 =  \{ (-\log\frac{j}{k}, \log \frac{X_{(j)}}{X_{(k)}}); 1 \le j
 \le k \}
            =  \{ (-\log t, \log \frac{X_{([kt])}}{X_{(k)}}) ; t \in \{\frac{1}{k}, \ldots,
           \frac{k-1}{k},1\}\},
$$
and
also write $\mathcal{T}$ as
$$
\mathcal{T} = \{(x, \frac{x}{\alpha}); x \ge 0\}
 = \{(-\log t, -\frac{1}{\alpha}\log t) ; 0 < t \le 1 \},
$$
where we put $x= - \log
     t$. We first show $\mathcal{S}'_n \stackrel{P}{\to} \mathcal{T}.$

Referring to Lemma \ref{lem:approxCurve}, set
$$x(t)=-\log t, \quad Y_n(t)=\log \frac{X_{(\lceil kt\rceil )}}{X_{(k)}},\quad
y(t)=-\frac{1}{\alpha} \log t,\qquad 0<t\leq 1.$$
>From \citet[page 82, equation (4.18)]{resnickbook:2006}, we have
$Y_n \stackrel{P}{\to} y,$ in $D_l(0,1]$, the left continuous
  functions on $(0,1]$ with finite right limits, metrized by the
    Skorohod metric.
Suppose $\{n''\}$ is a subsequence. There exists a
  further subsequence $\{n'\}\subset \{n''\}$ such that
$Y_{n'} \stackrel{a.s.}{\to} y$, in $D_l(0,1]$, and by
Lemma \ref{lem:approxCurve}, $\mathcal{S}'_{n'} \stackrel{a.s.}{\to}
\mathcal{T} $ in $\F$. Therefore $\mathcal{S}'_n \stackrel{P}{\to}
\mathcal{T}$, in $\F$, as $n \to \infty$.

 Now observe that with $\ba_n := (-\log{\frac{k}{n+1}},  \log
 X_{(k)}),$ we have
\begin{align*}
  \mathcal{S}_n  = & \{ (-\log{\frac{j}{n+1}},\log X_{(j)}) ; j
           =1,\ldots,k\} \\
           = & \{ (-\log{\frac{j}{k}} ,\log \frac{X_{(j)}}{X_{(k)}}  ) ; j =1,\ldots,k\} +  (-\log{\frac{k}{n+1}},  \log X_{(k)})\\
           = & \mathcal{S}_n^{'} + \ba_n.
  \end{align*}
Also,
$$\mathcal{T}_n  = \{(x,\frac{x}{\alpha}) ; x \ge 0 \} + (-\log
      \frac{k}{n+1},\log X_{(k)})\\
       =  \mathcal{T} + \ba_n.
$$
 Now, since   $d_\F(\mathcal{S}_n^{'},
 \mathcal{T}) \cinP  0$,  we get
$$
d_\F (\mathcal{S}_n,\mathcal{T}_n)  =  d_\F (\mathcal{S}_n^{'}+\ba_n,T+\ba_n)
           =  d_\F (\mathcal{S}_n^{'},T)
 \cinP 0,
$$ as required.
\end{proof}

   \section{Least squares line through a closed set} \label{sec:lsline}

  \subsection{Convergence of the least squares line}
  \label{subsec:conls}
The previous two sections gave results about the
   convergence of the QQ plot
to a straight line in the Fell topology of $\F_2.$
It is of interest to know whether some functional of  closed sets is
 continuous or not and, in particular, the slope of
   the least squares line through the points of QQ plot is one such
functional.  The slope of the least squares line is an estimator of
   scale for
   location/scale families and this leads to an estimate of the heavy
   tail index $\alpha$; see \citet{kratz:resnick:1996,
beirlant:vynckier:teugels:1996} and \citet[Section 4.6]{resnickbook:2006}.

Intuition suggests that when a sequence of finite sets converges to a
line, the slope of the least squares line should converge to the slope
of the limiting line. However there are subtleties which prevent this
from being true in general.  We need some restriction on the point
sets that converge, since  otherwise, a sequence of point sets which are
essentially linear
except for a vanishing bump, may converge to a line but the bump may
skew the least squares line sufficiently to prevent the slope from
converging; see Example \ref{eg:crit2} below.

The following Proposition provides a condition for the continuity
property to hold. First define the subclass $\Ffiniteorline \subset
\F_2$ to be the closed sets of $\F_2$ which are either sets of finite
cardinality or closed, bounded line segments.  These are the only
cases of compact sets where it is clear how to define a least squares
line. For $F \in \Ffiniteorline$, the functional $LS$ is defined in
the obvious way:
\begin{equation*}
    LS(F) = \text{slope of the least squares line through the closed set } F
    \end{equation*}

For the next proposition, we consider sets
$F_{n}:=\{ (x_i(n),y_i(n)): 1 \le i \le k_n\}$ of points
and write $\bar x_n=\sum_{j=1}^{k_n} x_j(n)/k_n$ and
$\bar y_n=\sum_{j=1}^{k_n} y_j(n)/k_n.$ Also, for a finite set $\mathcal{S}_n$,
$\#\mathcal{S}_n$ denotes the cardinality of $\mathcal{S}_n$.

      \begin{Proposition} \label{propn:lsline}
Suppose we have a sequence of sets $F_{n}:=\{ (x_i(n),y_i(n)): 1 \le i
\le k_n\} \in \K_2$, each consisting of $k_n$ points, which converge
to a bounded line segment $F \in \K_2$ with slope $m $ where $|m| <
\infty $, as $k_n\to\infty$.  Then
$$
LS(F_n) \to LS(F) = m
$$
provided the following condition holds:
\begin{align}
\exists  ~\delta > 0,  \text{ such that } p_{\delta}^n: =
\frac{ \#\Bigl(
  \{(\bar{x}_n-\delta,\bar{x}_n + \delta) \times
  (\bar{y}_n-\delta,\bar{y}_n + \delta)\} \bigcap F_n\Bigr)}{\# F_n} \to
p_{\delta } \in [0,1) .
\label{cond:notone}
      \end{align}
\end{Proposition}

 This Proposition gives a condition for the continuity of the slope
 functional $LS (\cdot)$
when $\{F_n n\ge 1\}$ and $ F$ are bounded sets in
 $\Ffiniteorline$.  The next example  shows the necessity of condition
 (\ref{cond:notone}), which prevents a  set of outlier points from
 skewing the slope of the least squares line.

 \begin{Example} \label{eg:crit2}
For $n\geq 1$, define the sets:
$$
  F_n  = \{(\frac{i}{n},0) , -n \le i \le n; (\frac{1}{n}(1+
\frac{j}{2^n}),\frac{1}{n}(1+ \frac{j}{2^n})), 0 \le j\le 2^n \}
\quad \text{ and }
\quad  F  =[-1,1] \times \{0\}.
$$
We develop features  about this example.
\begin{enumerate}
\item For the cardinality of $F_n$ we have
$$\#F_n=k_n=2^n +2n + 2.$$
\item  We have $F_n \to F$ in $\K_2$.
As before,
 denote the Hausdorff distance between two closed
sets in $\K_2$ by $D(\cdot,\cdot)$ and we have
$D(F_n, F) < {3}/{n} \to 0$ as $n \to \infty$.

\item Condition (\ref{cond:notone}) is {\it not\/} satisfied.
To see this pick any $n \ge 1$ and observe
$$
\bar{x}_n =  \bar{y}_n = \frac{3(2^n+1)}{2n(2^n+2n+2)}
=\frac{3(2^n+1)}{2n k_n}
\sim
\frac{3}{2n}.
$$
Fix $\delta > 0$. For all $n$ so large that $\delta >1/(2n)$ we have
$$ \frac{ \#\Bigl( \{(\bar{x}_n-\delta,\bar{x}_n + \delta) \times
  (\bar{y}_n-\delta,\bar{y}_n + \delta)\}  \bigcap F_n \Bigr)}{\# F_n} \ge
\frac{2^n+1}{2^n+2n+2}  \to 1,\quad (n\to\infty).$$
\end{enumerate}

\item Obviously for this example,
 $ m  = LS(F) = 0$. However,
if $m_n$ denotes the slope of the least squares line
through  $F_n$ then we show that $m_n \to 1 \neq 0=m$.
To see this, observe that conventional wisdom yields,
\begin{equation}\label{eqn:conventionalwisdom}
 m_n  = \frac
{\sum\limits_{(x_i(n),y_i(n)) \in F_n} (y_i(n) - \bar{y})(x_i(n) -
\bar{x})}{\sum\limits_{(x_i(n),y_i(n)) \in F_n} (x_i(n) -
 \bar{x})^2}.
\end{equation}
For the numerator we have,
\begin{align*}
 \sum\limits_{(x_i(n),y_i(n)) \in F_n} & (y_i(n) - \bar{y}_n)(x_i(n) -
  \bar{x}_n)  =
  \sum\limits_{(x_i(n),y_i(n)) \in F_n} y_i(n) x_i(n) - k_n \bar{y}_n
  \bar{x}_n  \\
 & =  \sum\limits_{j=0}^{2^n} \frac{1}{n^2}(1+ \frac{j}{2^n})^2
 - k_n \Bigl(  \frac{3(2^n+1)}{2n k_n} \Bigr)^2
= \frac{1}{n^2} \Biggl( \sum_{j=0}^{2^n}\Bigl(
1+\frac{2j}{2^n} +\frac{j^2}{2^{2n}} \Bigr)
-\frac{9}{4k_n} (2^n +1)^2
\Biggr)\\
&= \frac{1}{n^2} \Biggl(
2\cdot (2^n+1) +\frac{1}{2^{2n}} \sum_{j=0}^{2^n}j^2 -
\frac{9}{4k_n} (2^n +1)^2
\Biggr)\\
\intertext{and using the identity $\sum_{j=1}^N j^2=N(N+1)(N+\frac
  12)/3=N(N+1)(2N+1)/6,$ we get the above equal to}
&= \frac{1}{n^2} \Biggl(
2\cdot(2^n+1)+\frac{1}{2^{2n}} \frac{ 2^n(2^n+1)(2^n +\frac 12) }{3}
-\frac{9}{4k_n} (2^n +1)^2 \Biggr)\\
&=\frac{2^n+1}{n^2}\Bigl(
2+\frac{2^n +\frac 12 }{3\cdot 2^n}-\frac{9}{4k_n} (2^n +1)
\Bigr)
 \sim \frac{k_n}{12n^2}.
\end{align*}
For the denominator, we use the calculation already done for the numerator:
\begin{align*}
 \sum\limits_{(x_i(n),y_i(n)) \in F_n}& (x_i(n) - \bar{x}_n)^2  =
  \sum\limits_{(x_i(n),y_i(n)) \in F_n} x_i(n)^2 - k_n ( \bar{x}_n)^2
 \\
 & =  \sum\limits_{i=-n}^{n} (\frac{j}{n})^2 + \sum\limits_{j=0}^{2^n}
  \frac{1}{n^2}(1+ \frac{j}{2^n})^2
 -k_n \Bigl(\frac{3(2^n+1)}{2nk_n}\Bigr)^2 \\
&= \sum\limits_{i=-n}^{n} (\frac{j}{n})^2 +
  \sum\limits_{(x_i(n),y_i(n)) \in F_n} y_i(n) x_i(n) - k_n \bar{y}_n
  \bar{x}_n
\\
 & =  \frac{2n(n+1)(2n+1)}{6n^2} + \frac{k_n}{12n^2}
  +o(\frac{k_n}{12n^2})\\
&=O(n)+\frac{k_n}{12n^2} +o(\frac{k_n}{12n^2})\sim \frac{k_n}{12n^2}.
\end{align*}
Combining the asymptotic forms for numerator and denominator with
\eqref{eqn:conventionalwisdom} yields
$$
 m_n  \sim  \frac{k_n/12n^2}{k_n/12n^2} \sim 1,\quad (n\to\infty),$$
so $m_n \to 1 \ne 0 = m$, as claimed. \qed
 \end{Example}

\noindent {\it Proof of Proposition
      \ref{propn:lsline}\/}.
 For $(x_i(n),y_i(n)) \in F_n,$ we can write
\begin{equation} \label{eqn:lsyx}
       y_i(n) = m x_i(n) + z_i(n) \qquad 1\le i \le k_n
\end{equation}
We want to show that $m_n = LS (F_n) \to m =LS (F)$, as $n \to
\infty$.  Fix $\epsilon >0$. We will provide $N$ such that for $n> N$,
we have $|m_n - m | < \epsilon$.

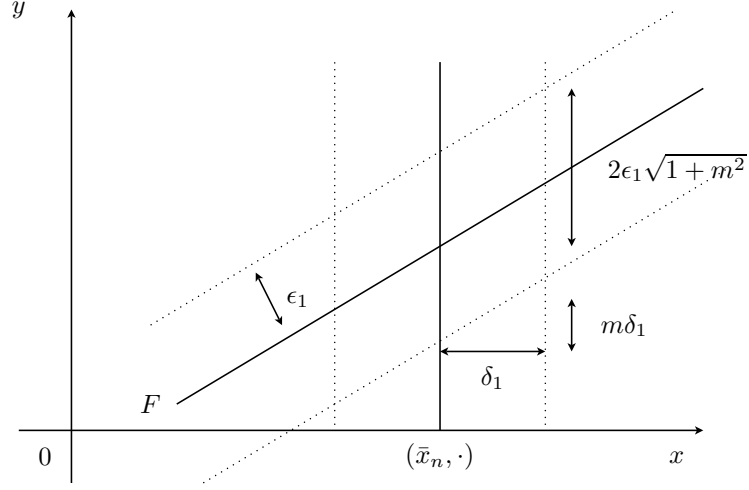
\begin{figure}[h]
\ifx\JPicScale\undefined\def\JPicScale{0.7}\fi
\psset{unit=\JPicScale mm}
\psset{linewidth=0.3,dotsep=1,hatchwidth=0.3,hatchsep=1.5,shadowsize=1,dimen=middle}
\psset{dotsize=0.7 2.5,dotscale=1 1,fillcolor=black}
\psset{arrowsize=1 2,arrowlength=1,arrowinset=0.25,tbarsize=0.7
5,bracketlength=0.15,rbracketlength=0.15}
\begin{pspicture}(0,0)(150,90)
\psline{<-}(20,90)(20,0) \psline{->}(10,10)(140,10)
\psline[linestyle=dotted](35,30)(135,90) \psline(40,15)(140,75)
\psline[linestyle=dotted](45,0)(145,60)
\rput{0}(150,70){\psellipse[](0,0)(0,0)} \psline(90,80)(90,10)
\psline[linestyle=dotted](70,80)(70,10)
\psline[linestyle=dotted](110,80)(110,10)
\rput(90,5){$(\bar{x}_n,\cdot)$} \psline{<->}(90,25)(110,25)
\psline{<->}(115,35)(115,25) \psline{<->}(115,75)(115,45)
\rput(100,20){$\delta_1$} \rput(125,30){$m\delta_1$}
\rput(135,60){$2\epsilon_1\sqrt{1+m^2}$} \rput(100,20){}
\psline{<->}(55,40)(60,30) \rput(63,35){$\epsilon_1$}
\rput(115,30){} \rput(35,15){$F$} \rput(135,5){$x$}
\rput(10,90){$y$} \rput(15,5){0}
\end{pspicture}
\caption{The geometry of a neighborhood of the line; support for
  \eqref{eqn:delta}. }
\label{fig:two}
\end{figure}

 First of all,  condition (\ref{cond:notone}) allows us to
 fix $\delta > 0       $ such that
$$
p_\delta^n:= p_n = \frac{ \# \{(\bar{x}_n-\delta,\bar{x}_n +
      \delta) \times (\bar{y}_n-\delta,\bar{y}_n + \delta)\} \bigcap
      F_n}{\# F_n} \to p < 1 .
$$
Choose $N_1$ such that for $n > N_1$, we have $p_n < \frac{1+p}{2}$
or equivalently that $1-p_n > \frac{1-p}{2}$. For $\eta > 0$ and $F
\in \K_2$, recall the definition of the $\eta$-swelling of $F$:
\begin{equation}\label{def:Fepsil}
     F^{\eta} = \{x: d(x,y)< \eta \mathrm{~~for~~ some~~} y \in F \}.
\end{equation}
Since $D(F_n, F )\to 0$ in $\K_2$,
we can  choose $N_2$ such that for all
$n > N_2$ we have $F_n \subset
      F^{\epsilon_1}$ where
$$
 \epsilon_1:= \frac{2 \delta\epsilon (1-p)}{4\sqrt{1+m^2}(2+2m+
 \epsilon(1-p))}
=\delta_1\epsilon\frac{(1-p)}{4\sqrt{1+m^2}}
$$
and we have set
$$
\delta_1:= \frac{\delta}{1+m +\frac 12 \epsilon (1-p)} <\delta.
$$
The choice of $\delta_1$ is designed to ensure that if for some
$(x_i(n),y_i(n))$,
 we have $|x_i(n)-\bar{x}_n| < \delta_1$, then
$$
(x_i(n),y_i(n)) \in
(\bar{x}_n-\delta, \bar{x}_n +\delta) \times (\bar{y}_n-\delta,
 \bar{y}_n +\delta).$$
This follows  because
\begin{align}
|x_i(n)-\bar{x}_n|\vee & |y_i(n)-\bar{y}_n|
 < \delta_1 + m \delta_1 + 2\epsilon_1\sqrt{1+m^2} .\nonumber \\
\intertext{See  Figure \ref{fig:two};
 from the definition of $\epsilon_1$ we have this equal to}
  & = \delta_1 + m \delta_1 +
 2\frac{\delta_1\epsilon(1-p)}{4\sqrt{1+m^2}}\sqrt{1+m^2}
 =\delta_1(1 + m +
 \frac{\epsilon(1-p)}{2}) =\delta. \label{eqn:delta}
        \end{align}

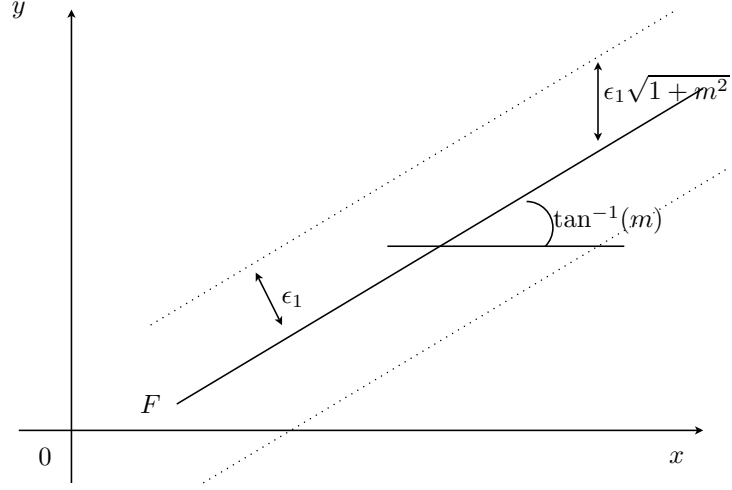
\begin{figure}[t]
\ifx\JPicScale\undefined\def\JPicScale{0.7}\fi
\psset{unit=\JPicScale mm}
\psset{linewidth=0.3,dotsep=1,hatchwidth=0.3,hatchsep=1.5,shadowsize=1,dimen=middle}
\psset{dotsize=0.7 2.5,dotscale=1 1,fillcolor=black}
\psset{arrowsize=1 2,arrowlength=1,arrowinset=0.25,tbarsize=0.7
5,bracketlength=0.15,rbracketlength=0.15}
\begin{pspicture}(0,0)(150,90)
\psline{<-}(20,90)(20,0) \psline{->}(10,10)(140,10)
\psline[linestyle=dotted](35,30)(135,90) \psline(40,15)(140,75)
\psline[linestyle=dotted](45,0)(145,60)
\rput{0}(150,70){\psellipse[](0,0)(0,0)}
\psline{<->}(120,80)(120,65) \rput(100,20){}
\psline{<->}(55,40)(60,30) \rput(62,35){$\epsilon_1$}
\rput(115,30){} \rput(35,15){$F$} \rput(135,5){$x$}
\rput(10,90){$y$} \rput(15,5){0} \rput(133,75){$\epsilon_1\sqrt{1+
m^2}$} \rput(90,60){} \psline(80,45)(125,45)
\rput{0}(106.46,48.54){\psellipticarc[](0,0)(5,-5){-90}{45}}
\rput(122,50){$\tan^{-1}( m)$}
\end{pspicture}
\caption{The geometry of the neighborhood of the line $F$.}
\label{fig:one}
\end{figure}

  Let $ N = N_1 \vee N_2 $ and restrict attention to  $n > N$.
 Since $F_n\subset F^{\epsilon_1}$, we have for all $1 \leq i \leq k_n$
 that
  $(x_i(n),y_i(n)) \in F^{\epsilon_1}. $  By convexity of
$F^{\epsilon_1}$,  $(\bar{x},\bar{y}) \in F^{\epsilon_1}
         $.
Therefore, referring to Figure \ref{fig:one}, we have
\begin{align}
|z_i(n) - \bar{z}_n| \le & |y_i(n) - m x_i(n)| + |\bar{y}_n- m \bar{x}_n|
                       \nonumber \\
\le & \epsilon_1 \sqrt{1+m^2} + \epsilon_1
                       \sqrt{1+m^2} = 2\epsilon_1 \sqrt{1+m^2}.\label{eqn:ub}
           \end{align}
Using the representation (\ref{eqn:lsyx}) we get,
\begin{equation}\label{eqn:repm}
 m_n  =  \frac{\sum\limits_{i=1}^{k_n}
 (y_i(n)-\bar{y}_n)(x_i(n)-\bar{x}_n)}{\sum\limits_{i=1}^{k_n}
                    (x_i(n)-\bar{x}_n)^2}
=  m + \frac{ \sum\limits_{i=1}^{k_n}(z_i(n) - \bar{z}_n)(x_i(n)-\bar{x}_n)}{\sum\limits_{i=1}^{k_n}
                     (x_i(n)-\bar{x}_n)^2}.
\end{equation}
 Therefore,
\begin{align*}
|m_n - m | & =
\Biggl|\frac{
        \sum\limits_{i=1}^{k_n}(z_i(n) -
        \bar{z}_n)(x_i(n)-\bar{x}_n)}{\sum\limits_{i=1}^{k_n}
                     (x_i(n)-\bar{x}_n)^2}\Biggr|
 \le  \frac{ \sum\limits_{i=1}^{k_n} |z_i(n) -
   \bar{z}_n||x_i(n)-\bar{x}_n|}{\sum\limits_{i=1}^{k_n}
                     (x_i(n)-\bar{x}_n)^2}
                       \le 2 \epsilon_1 \sqrt{1+m^2} \frac{ \sum\limits_{i=1}^{k_n} |x_i(n)-\bar{x}_n|}{\sum\limits_{i=1}^{k_n}
                     (x_i(n)-\bar{x}_n)^2} \\
       \end{align*}
where the last inequality  follows from \eqref{eqn:ub}.

For convenience, define the following  notation:
         \begin{align*}
|S(x)|_{< \rho} &:= \sum\limits_{|x_i(n)-\bar{x}_n| < \rho}
          |x_i(n)-\bar{x}_n|, &
|S(x)|_{\ge \rho} &:= \sum\limits_{|x_i(n)-\bar{x}_n| \ge
        \rho} |x_i(n)-\bar{x}_n| , \\
S^2(x)_{< \rho} &:= \sum\limits_{|x_i(n)-\bar{x}_n| < \rho}
   (x_i(n)-\bar{x}_n)^2,
   &
S^2(x)_{\ge \rho} &:= \sum\limits_{|x_i(n)-\bar{x}_n| \ge \rho}
(x_i(n)-\bar{x}_n)^2 ,\\
B\bigl( (x,y),\delta\bigr)&:= ({x}-\delta,{x} +
      \delta) \times ({y}-\delta,{y} + \delta). & {}
         \end{align*}

Therefore
\begin{align*}
\frac{ \sum\limits_{i=1}^{k_n} |x_i(n)-\bar{x}_n|}{\sum\limits_{i=1}^{k_n}
  (x_i(n)-\bar{x}_n)^2}
& =  \frac{|S(x)|_{<\delta_1}+ |S(x)|_{\ge
  \delta_1}}{S^2(x)_{<\delta_1}+ S^2(x)_{\ge \delta_1}}
=  \frac{(|S(x)|_{<\delta_1}+ |S(x)|_{\ge \delta_1})/S^2(x)_{\ge
  \delta_1}}{(S^2(x)_{<\delta_1}/S^2(x)_{\ge\delta_1}+1)}\\
& \le  \frac{(|S(x)|_{<\delta_1}+ |S(x)|_{\ge\delta_1})}{S^2(x)_{\ge
      \delta_1}}
\le  \frac{1}{\delta_1}
  \biggl(\frac{|S(x/\delta_1)|_{<1}}{S^2(x/\delta_1)_{\ge 1}} +
  1\biggr)\\[2mm]
 &\le  \frac{1}{\delta_1}
   \biggl(\frac{\# \{(x_i(n),y_i(n))\in F_n: |x_i(n)-\bar{x}_n| <
  \delta_1\}} {\# \{(x_i(n),y_i(n))\in F_n: |x_i(n)-\bar{x}_n| \ge
  \delta_1\}} +  1\biggl)\\[2mm]
& \le  \frac{1}{\delta_1} \biggl(\frac{\# \{(x_i(n),y_i(n))\in F_n:
  (x_i(n),y_i(n)) \in B((\bar{x}_n,\bar{y}_n),\delta)\}} {\#
  \{(x_i(n),y_i(n))\in F_n: (x_i(n),y_i(n)) \notin
  B((\bar{x}_n,\bar{y}_n),\delta)\}} +  1\biggl)\\[2mm]
\intertext{The choice of $\delta_1$ justifies the previous
  step by \eqref{eqn:delta}. The previous expression is bounded by}
& \le  \frac{1}{\delta_1} (\frac{p_n}{1-p_n} + 1)
 \le  \frac{1}{\delta_1} (\frac{1+p}{1-p} + 1)
 =  \frac{2}{\delta_1(1-p)},
\end{align*}
and we recall $p<1$.

Consequently
  \begin{align*}
|m_n - m |  = & 2 \epsilon_1 \sqrt{1+m^2} \frac{ \sum\limits_{i=1}^{k_n}
|x_i(n)-\bar{x}_n|}{\sum\limits_{i=1}^{k_n} (x_i(n)-\bar{x}_n)^2}
    \le  2 \epsilon_1 \sqrt{1+m^2} \times \frac{2}{\delta_1(1-p)}\\
    = & 2 \epsilon\delta_1\frac{1-p}{4\sqrt{1+m^2}} \sqrt{1+m^2} \times \frac{2}{\delta_1(1-p)}= \epsilon
       \end{align*}
This completes the proof that  $m_n \to m $ under condition
\eqref{cond:notone}.  \qed

\medskip

\begin{Corollary} \label{cor:means}
If $\bar{x}_n \to \mu_x < \infty$ and $\bar{y}_n \to \mu_y <
\infty$, as $n \to \infty$, then Proposition \ref{propn:lsline}
 holds if we replace $(\bar{x}_n,\bar{y}_n)$ in (\ref{cond:notone})
by $(\mu_x,\mu_y)$.
\end{Corollary}

\begin{proof}
 In place of condition (\ref{cond:notone}) we are assuming
 \begin{align}
      \exists  ~\delta > 0  \text{ such that } p_{\delta}^n = \frac{
      \# \{(\mu_x-\delta,\mu_x + \delta) \times (\mu_y-\delta,\mu_y +
      \delta)\} \bigcap F_n}{\# F_n} \to p_{\delta } \in [0,1).
 \label{eqn:modCondit}
     \end{align}
Let us fix $\delta > 0$ such that
$$ p_n^{*} := \frac{ \# \{(\mu_x-2\delta,\mu_x + 2\delta) \times
(\mu_y-2\delta,\mu_y + 2\delta)\} \bigcap F_n}{\# F_n} \to p \in [0,1).
$$
Since $\bar{x}_n \to \mu_x < \infty$ and $\bar{y}_n \to \mu_y <
\infty$, there exists $N^{*}$ such that $n > N^{*}$ implies that
$(\bar{x}_n,\bar{y}_n) \in (\mu_x-\delta,\mu_x + \delta) \times
(\mu_y-\delta,\mu_y + \delta)$. Hence for $n > N^{*}$
\begin{align*}
 p_n & : = \frac{ \# \{(\bar{x}_n-\delta,\bar{x}_n +
\delta) \times (\bar{y}_n-\delta,\bar{y}_n +
\delta)\} \bigcap F_n}{\# F_n}\\
& \le \frac{ \# \{(\mu_x-2\delta,\mu_x + 2\delta) \times
(\mu_y-2\delta,\mu_y + 2\delta)\} \bigcap F_n}{\# F_n}\\
 & = p_n^{*} \to p \in [0,1).
\end{align*}
 Now choose $N_1 \ge N^{*}$ such that for all $n > N_1$, we have $p_n <
      \frac{1+p}{2}$. This also means that $1-p_n  > \frac{1-p}{2}$.

\noindent The rest of the proof is  the same as that of Proposition
\ref{propn:lsline}.
\end{proof}

\section{Slope of the LS line as a tail index estimator}
For heavy tailed distributions, the slope of the least squares line
 through the QQ plot made by the upper $k_n$ largest order statistics
 is a consistent estimator of $1/\alpha$.
 See \citet{kratz:resnick:1996,
beirlant:vynckier:teugels:1996} and \citet[Section
 4.6]{resnickbook:2006}. We connect the ideas of the previous section
 with this result.

\begin{Proposition}\label{prop:consistentEstr}
Consider non-negative random variables $ X_1,X_2,\ldots,X_n$ which are
iid with common distribution
$F$ where $\bar F
\in RV_{-\alpha}$ and $ X_{(1)} \ge X_{(2)} \ge \ldots \ge X_{(n)}$
are the order statistics in decreasing order.
The sets  $\mathcal{S}_n$ and $\mathcal{T}_n$ were
 defined in Proposition \ref{propn:conrv} where we proved
$d_\F (\mathcal{S}_n,\mathcal{T}_n)\stackrel{P}{\to} 0, $ assuming
 $k=k(n)\to\infty$ and $k/n \to 0$ as $n \to \infty$.
For  convenience we   defined
$\snp = \mathcal{S}_n + a_n$ and $\mathcal{T}=
\mathcal{T}_n + a_n$ where $a_n$ was a random point.
Write
$$
\mathcal{S}_n^{'}  = \{ (-\log{\frac{j}{k_n}},
\log \frac{X_{(j)}}{X_{(k)}} ; j=1,\ldots,k_n \}
  = \{ (x_j(n),y_j(n)) ; j
=1,\ldots,k_n \} \;\text{( say)} \text{ and }
\mathcal{T}  = \{(x,\frac{x}{\alpha}) ; x \ge 0 \}.$$
Then,
 \begin{align} \label{eqn:alconrv}
LS(\mathcal{S}_n^{'}) = {LS}(\mathcal{S}_n) \cinP \frac{1}{\alpha} =
 {LS}(\mathcal{T}_n)=
 {LS}(\mathcal{T}),
\end{align}
as $k:=k_n \to \infty$ and $k_n/n \to 0$ as $n \to
\infty$.
\end{Proposition}

The result is believable  based on the fact that $d_\F
(\mathcal{S}_n,\mathcal{T}_n)\stackrel{P}{\to} 0.$
 However, since neither $\mathcal{T}_n$ nor $\mathcal{T}$ are $\K_2$
  sets, some sort of truncation to compact regions of $\mathbb{R}^2$
  is necessary in order to capitalize on Proposition
  \ref{propn:lsline}. For some integer $M  > 2$, define
$$ K_M  = [0,M] \times [0, \frac{2M}{\alpha}],
$$
and let
$${\mathcal{S}_n^{'}}^M=\mathcal{S}_n^{'}\cap K_M \quad
  \text{and}\quad \mathcal{T}^M=\mathcal{T}\cap K_M.$$

\begin{proof}
Some preliminary observations. Clearly,
$
LS(\mathcal{S}_n) =LS(\snp + a_n) = LS(\snp)
$
and with $x_j(n), y_j(n)$ defined in the statement of the
Proposition,
$$
  LS(\snp)  =
 \frac{\bar{S}_{XY}-\bar{S}_{X}\bar{S}_{Y}}{\bar{S}_{XX}-(\bar{S}_{X})^2},
$$
where,  as usual,
\begin{align*}
\bar{S}_X & = \frac{1}{k_n}\sum\limits_{(x_j(n),y_j(n)) \in \snp }
x_j(n) , & \bar{S}_Y & =
\frac{1}{k_n}\sum\limits_{(x_j(n),y_j(n)) \in \snp } y_j(n), \\
\bar{S}_{XY} & = \frac{1}{k_n}\sum\limits_{(x_j(n),y_j(n)) \in
\snp } x_j(n) y_j(n) , & \bar{S}_{XX} & =
\frac{1}{k_n}\sum\limits_{(x_j(n),x_j(n)) \in \snp } (x_j(n))^2 .
\end{align*}
We need similar quantities $\bar S_X^M, \bar S_Y^M, \bar S_{XY}^M$
corresponding to averages of points restricted to  $K_M$, so for instance
$$\bar{S}_X^M  = \frac{1}{k^M}\sum\limits_{(x_j(n),y_j(n)) \in \SM }
x_j(n) $$
and $k^M=\#\SM$. A simple calculation given in \citet[page
  109]{resnickbook:2006} yields as
$k \to \infty$,
\begin{equation}\label{eqn:convergences}
\bar{S}_X   = \frac{1}{k} \sum\limits_{i=1}^k (-\log\frac{i}{k})
\sim \int_0^1 (-\log x) dx =1,\quad
\bar{S}_{XX}  = \frac{1}{k} \sum\limits_{i=1}^k
(-\log\frac{i}{k})^2 \sim  \int_0^1 (-\log x)^2 dx
  =2,
\end{equation}
while for $\bar S_Y$ we have
\begin{equation}\label{eqn:Hillconv}
\bar{S}_Y=
 \frac{1}{k} \sum\limits_{i=1}^k
(-\log\frac{X_{(i)}}{X_{(k)}}) \cinP \frac{1}{\alpha}
\end{equation}
since $\bar{S}_Y$ is the Hill estimator and is consistent for
$1/\alpha$
\citep{resnickbook:2006, csorgo:deheuvels:mason:1985, mason:1982,
  mason:turova:1994}.

We need the corresponding limits for $\bar S_X^M, \bar S_{XX}^M, \bar
S_Y^M$. These calculations and subsequent calculations are simplified
by the following facts:
\begin{enumerate}
\item The ratios of order statistics process converges, as $k \to
  \infty$, $k/n\to 0$,
\begin{equation}\label{eqn:ratioConv}
\frac{X_{(\lceil k t\rceil)}}{X_{(k)}} \cinP
t^{-1/\alpha},
\end{equation}
in $D_l(0,\infty]$ \citep[page 82]{resnickbook:2006}.
\item Define the random measure
$$\hat \nu_n (\cdot) =\frac 1k \sum_{i=1}^n \epsilon_{X_{(i)}/X_{(k)}
} (\cdot)$$
on $(0,\infty]$, which puts mass $1/k$ at the points
  $\{X_{(i)}/X_{(k)}, 1 \leq i \leq n \}$. Then
\begin{equation}\label{eqn:nuHatConv}
\hat \nu_n \cinP \nu_\alpha,
\end{equation} in the space of Radon measures on
    $(0,\infty]$, where $\nu_\alpha (x,\infty]=x^{-\alpha},\,x>0$
    \citep[page 83]{resnickbook:2006}.
\item The number of points $k^M $ in $\SM$ satisfies, as $n \to
  \infty$,  $k\to
  \infty$, $k/n\to 0$,
\begin{equation}\label{eqn:kMform}
k^M /k \cinP 1-e^{-M}.
\end{equation}
To see this, observe
\begin{align*}
k^M/k =&\frac 1k \#\{j\leq k: k\geq j\geq ke^{-M} \text{ and }
\frac{X_{(j)}}{X_{(k)}}  \leq e^{2M/\alpha} \}\\
=&\frac 1k \# \{j\leq k : 1\leq \frac{X_{(j)}}{X_{(k)}}
\leq \frac{X_{(\lceil ke^{-M}  \rceil)}}{X_{(k)}} \wedge
e^{2M/\alpha} \}\\
=& \hat \nu_n \Bigl( 1, \frac{X_{(\lceil ke^{-M}  \rceil)}}{X_{(k)}} \wedge
e^{2M/\alpha} \Bigr]\\
\cinP & \, 1-\Bigl( (e^{-M})^{-1/\alpha} \wedge
e^{2M/\alpha}\Bigr)^{-\alpha}=1-e^{-M}.
\end{align*}
\end{enumerate}

We continue using these three facts.
For $\bar S_X^M$ we have
$$ \bar S_x^M=\frac{1}{k^M} \sum_{\bigl( x_i(n),y_i(n)\bigr) \in \SM}
x_i(n)
=\frac{1}{k^M}
\sum_{\substack{j:k\geq j\geq ke^{-M}\\
0<\log X_{(j)}/X_{(k)} \leq 2M/\alpha }}
-\log \frac jk.
$$
Set
\begin{align*}
\Bigl(\bar S_X^M \Bigr)^* :=&\frac{1}{k^M}
\sum_{j:k\geq j\geq ke^{-M}}-\log \frac jk
=\frac{k}{k^M} \frac 1k \sum_{j:k\geq j\geq ke^{-M}}-\log \frac jk\\
\sim & \frac{1}{1-e^{-M}} \int_{e^{-M}}^1 -\log x \,dx =
\frac{1}{1-e^{-M}} \int_{0}^M  y e^{-y} dy\\
=:& 1+\epsilon_X (M),
\end{align*}
where $\epsilon_X(M) \to 0$ as $M \to \infty$.
Also, $\bar S_X^M$ and
$\Bigl(\bar S_X^M \Bigr)^* $  are close asymptotically since
\begin{align*}
P[\bar S_X^M \neq
\Bigl(\bar S_X^M \Bigr)^* ]
=&P\Bigl\{ \bigcup_{k\geq j\geq k^{-M}} [\log \frac{X_{(j)}}{X_{(k)}}
  >2M/\alpha ]\Bigr\}\\
=&P[\log
\frac{X_{(\lceil ke^{-M}  \rceil)}}{X_{(k)}} >2M/\alpha]
\to 0 ,
\end{align*}
since
$$\frac{X_{(\lceil ke^{-M}  \rceil)}}{X_{(k)}}\cinP
e^{M/\alpha}<e^{2M/\alpha}.$$
We conclude
\begin{equation}\label{eqn:sxM}
\bar S_X^M \cinP 1+\epsilon_X(M):=\mu_X^M,
\end{equation}
with $\epsilon_X(M) \to 0$ as $M\to\infty$,
and in a similar way we can derive that
\begin{equation}\label{eqn:sxxM}
\bar S_{XX}^M \cinP 2+\epsilon_{XX}(M),
\end{equation}
where $\epsilon_{XX}(M) \to 0$ as $M\to\infty$.
For $\bar S_Y^M$ we have
\begin{align*}
\bar S_Y^M=&\frac{1}{k^M}
\sum_{
\substack{
j:k\geq j\geq ke^{-M}\\
0<\log X_{(j)}/X_{(k)} \leq 2\alpha^{-1} M
}
} \log \frac{X_{(j)}}{X_{(k)}}\\
=&\frac{1}{k^M}
\sum_{
j:0<\log X_{(j)}/X_{(k)}
 \leq 2\alpha^{-1}M     \wedge
\log X_{(\lceil ke^{-M} \rceil)}/X_{(j)}        }
\log \frac{X_{(j)}}{X_{(k)}}\\
=&\frac{k}{k^M} \int_1^{2\alpha^{-1} M\wedge
\log X_{(\lceil ke^{-M} \rceil)}/X_{(j)} } \log y \, \hat \nu_n (dy)\\
\cinP & \frac{1}{1-e^{-M}} \int_1^{2\alpha^{-1}M \wedge \alpha^{-1}M} \log y
\,\nu_\alpha (dy)\\
=& \frac{1}{1-e^{-M}}\int_0^{M/\alpha} se^{-\alpha s}ds =:\mu_Y^M,
\end{align*}
where $\mu_Y^M \to \frac{1}{\alpha} $ as $M \to \infty$. We conclude
\begin{equation}\label{eqn:barSYM}
\bar S_Y^M \cinP \mu_Y^M .
\end{equation}

To prove \eqref{eqn:alconrv}, we follow the following outline of steps.
\begin{itemize}
\item  Step 1: Prove $\SM \cinP \mathcal{T}^M$.
\item {Step 2:} Verify that Corollary \ref{cor:means} is applicable by
  showing that the analogue of \eqref{cond:notone} holds. This permits
  the conclusion that
$$LS(\SM) \cinP 1/\alpha. $$
Coupled with \eqref{eqn:sxM}, \eqref{eqn:sxxM} and \eqref{eqn:barSYM},
this yields
\begin{equation}\label{eqn:SxyM}
\bar S_{XY}^M=\frac{2}{\alpha} + \epsilon_{XY}(M)+o_p(1),
\end{equation}
where $\lim_{M\to \infty} \epsilon_{XY}(M)=0$ and $o_p(1) \cinP 0$ as
$n \to \infty$.
\item Step 3: Compare $\bar S_{XY}$ and $\bar S_{XY}^M$ and Check that
\begin{equation}\label{eqn:Compare}
\lim_{M\to \infty}\limsup_{n\to \infty} P[|\bar S_{XY}^M -\bar
  S_{XY}|>\eta] =0,\quad \forall \eta>0.
\end{equation}
This gives $\bar S_{XY}\cinP 2/\alpha$ which
coupled with \eqref{eqn:convergences} and
\eqref{eqn:Hillconv} implies \eqref{eqn:alconrv}.
\end{itemize}

We may check Step 1 using a very minor modification of Lemma
\ref{lem:approxCurve},  following the pattern of proof used for
Proposition \ref{propn:conrv}. For Step 2, the challenge is to verify
condition \eqref{cond:notone} holds and we defer this to the end of the
proof. Thus we turn to Step 3.

First of all, we observe that $\bar S_{XY}^M$ and $\bar S_{XY}$
average, respectively $k^M $ and  $k$ terms  but there is no need to
differentiate: For any $\eta>0$,
\begin{align*}
P[\Bigl|\frac{1}{k^M}\sum_{(x_j(n),y_j(n)) \in \SM} x_i(n)y_i(n)
-&\frac{1}{k}\sum_{(x_j(n),y_j(n)) \in \SM} x_i(n)y_i(n)\Bigr|>\eta]\\
=&P[\Bigl|\frac{1}{k^M}-\frac 1k \Bigr|
\sum_{(x_j(n),y_j(n)) \in \SM} x_i(n)y_i(n)>\eta]\\
\intertext{and dividing the sum by $k^M$ yields}
=&P\Bigl[\bar S_{XY}^M \Bigl|1-\frac{k^M}{k}\Bigr|>\eta\Bigr].
\end{align*}
Since $\bar S_{XY}^M$ is convergent in probability, it is
stochastically bounded and since, as $n\to\infty$,
$$\Bigl|1-\frac{k^M}{k}\Bigr|\;\stackrel{P}{\to}\; 1-(1-e^{-M})=e^{-M}
\; \stackrel{M\to\infty}{\to}\; 0,$$ we conclude
\begin{equation}\label{eqn:nibble1}
\lim_{M\to \infty}\limsup_{n\to \infty}
P\Bigl[\Bigl|\frac{1}{k^M}\sum_{(x_j(n),y_j(n)) \in \SM} x_i(n)y_i(n)
-\frac{1}{k}\sum_{(x_j(n),y_j(n)) \in \SM} x_i(n)y_i(n)\Bigr|>\eta\Bigr]
=0.
\end{equation}

Next observe for $\eta >0$,
\begin{align}
P[\Bigl|
 \frac{1}{k}\sum\limits_{(x_j(n),y_j(n)) \in \SM} x_j(n)
 y_j(n)
-&\frac 1k \sum_{k\geq j \geq ke^{-M}} x_j(n)y_j(n)\Bigr|>\eta]
\leq  P\{\bigcup_{k\geq j\geq
 ke^{-M}}[\frac{X_{(j)}}{X_{(k)}}>e^{2M/\alpha}]\} \nonumber \\
\leq & P[\frac{X_{(\lceil ke^{-M}  \rceil)}}{X_{(k)}}>e^{2M/\alpha}] \to
 0,\quad  (n\to\infty). \label{eqn:nibble2}
\end{align}
Note that by the Cauchy-Schwartz inequality,
$$\bigl(|\bar S_{XY} -\frac 1k \sum_{k\geq j\geq ke^{-M}} x_j(n)y_j(n)|\bigr)^2 \leq \bigl(\frac 1k \sum_{1\leq j
  \leq ke^{-M}} x_j(n)y_j(n)\bigr)^2
\leq
\frac 1k \sum_{1\leq j
  \leq ke^{-M}} x_j(n)^2 \; \cdot\;
\frac 1k \sum_{1\leq j
  \leq ke^{-M}} y_j(n)^2.
$$ Furthermore
\begin{align}
\frac 1k \sum_{1\leq j
  \leq ke^{-M}} y_j(n)^2
=&\int^\infty_{\log X_{(\lceil ke^{-M}   \rceil)}/X_{(k)}}    ( \log y
  )^2 \hat \nu_n (dy) \nonumber \\
\intertext{and using \eqref{eqn:ratioConv}, we have for some $c>0$,
all  large $n$ and some $M$ that the above is bounded by}
{}& \int^\infty_{cM}    ( \log y
  )^2 \hat \nu_n (dy) +o_p(1). \label{eqn:nibble3}
\end{align}
Assessing \eqref{eqn:nibble1}, \eqref{eqn:nibble2} and
\eqref{eqn:nibble3}, we see that
\eqref{eqn:Compare} will be proved if we show
\begin{equation}\label{eqn:reduction37}
\lim_{M\to\infty}\limsup_{n\to\infty}P\Bigl[ \int^\infty_{M}    ( \log y
  )^2 \hat \nu_n (dy) >\eta \Bigr]=0,\quad (\forall \eta >0).
\end{equation}
This treatment is similar to the stochastic version of Karamata's
theorem
(\citet{feigin:resnick:1997}, \citet[page 207]{resnickbook:2006}. For
  $0<\zeta <1\wedge \alpha $ and large $M$, the integrand $(\log y)^2 $
is dominated by $y^\zeta$. Bound the integral by
$$\int_M^\infty \hat \nu_n (y,\infty] \,\zeta y^{\zeta -1} dy
+M^\zeta \hat \nu_n (M,\infty].$$
If we let first $n\to\infty$ and then $M\to\infty$, for the second
piece we have
$$M^\zeta \hat \nu_n (M,\infty] \cinP M^\zeta \nu_\alpha (M,\infty]
    =M^{\zeta - \alpha} \to 0.$$

Now we deal with the integral. Set  $b(t)=(1/(1-F))^\leftarrow (t)$ so
that
$X_{(k)}/b(n/k) \cinP 1$ \citep[page 81]{resnickbook:2006}. For $\gamma >0$,
\begin{align*}
P[\int_M^\infty \hat \nu_n (y,\infty]\, \zeta y^{\zeta -1} dy> \eta]
=&P[\int_M^\infty \hat \nu_n (y,\infty] \,\zeta y^{\zeta -1} dy> \eta,
1-\gamma <X_{(k)}/b(n/k)<1+\gamma] + o(1)\\
\leq &
P[\int_M^\infty \frac 1k \sum_{i=1}^n \epsilon_{X_i/b(n/k)}
  ((1-\gamma) y,\infty] \zeta y^{\zeta -1} dy> \eta] +o(1).\\
\intertext{Ignore the term $o(1)$.  Markov's inequality gives a bound}
\leq & (const)\int_M^\infty  E\Bigl(\frac 1k \sum_{i=1}^n P[X_i\geq
  b(n/k)(1-\gamma) y] \Bigr) \zeta y^{\zeta -1} dy\\
= & (const)\int_M^\infty  \frac nk \bar F(
  b(n/k)(1-\gamma) y] ) \zeta y^{\zeta -1} dy.
\intertext{and applying
Karamata's theorem \citep{resnickbook:2006,
bingham:goldie:teugels:1987, geluk:dehaan:1987, dehaan:1970}, we have
as $n \to\infty$ that this converges to}
= & (const)\int_M^\infty  \bigl((1-\gamma) y )^{-\alpha} \zeta
y^{\zeta -1} dy \;\stackrel{M\to\infty}{\to} \; 0,
\end{align*}
as required. This finishes Step 3 and completes the proof modulo the
verification that \eqref{eqn:modCondit} can be proven for this problem.

The remaining task of checking \eqref{eqn:modCondit} proceeds as
follows. Recall $\mu_X^M$ and $\mu_Y^M$ from \eqref{eqn:sxM} and
\eqref{eqn:barSYM}. Fix $M$. Then for $p_\delta^n$ in
\eqref{eqn:modCondit}, we have
$$
\frac{1}{k^M} \#\{j:\mu_X^M-\delta <-\log \frac jk<\mu_X^M+\delta,
0<-\log \frac jk \leq M; \mu_Y^M-\delta<\log \frac{X_{(j)}}{X_{(k)}}
<\mu_Y^M+\delta, 0\leq \log \frac{X_{(j)}}{X_{(k)}}\leq \frac{ 2M}{\alpha}
\}.
$$
Since $\mu_X^M \approx 1$ and $\mu_Y^M\approx 1/\alpha$, we get for
large $M$
\begin{align*}
p_\delta^n:=&
\frac{1}{k^M} \#\Bigl\{j:\mu_X^M-\delta <-\log \frac jk<\mu_X^M+\delta;\;
 \mu_Y^M-\delta<\log \frac{X_{(j)}}{X_{(k)}}
<\mu_Y^M+\delta \Bigr\}\\
=&\frac{1}{k^M} \#\Bigl\{j: \frac{X_{(\lceil k\exp\{-(\mu_X^M-\delta )  \}
    \rceil)}}{X_{(k)}} \vee e^{\mu_Y^M -\delta}
< \frac{X_{(j)}}{X_{(k)}}
< \frac{X_{(\lceil k\exp\{-(\mu_X^M+\delta )  \}
    \rceil)}}{X_{(k)}} \wedge e^{\mu_Y^M+\delta}\Bigr\}\\
=&\frac{k}{k^M} \,\hat \nu_n \Biggl(
\frac{X_{(\lceil k\exp\{-(\mu_X^M-\delta )  \}
    \rceil)}}{X_{(k)}} \vee e^{\mu_Y^M -\delta},
 \frac{X_{(\lceil k\exp\{-(\mu_X^M+\delta )  \}
    \rceil)}}{X_{(k)}} \wedge e^{\mu_Y^M+\delta}
\Biggr).
\end{align*}
Apply \eqref{eqn:ratioConv} and \eqref{eqn:nuHatConv} and we find
\begin{align*}
p_\delta^n \cinP & \frac{1}{1-e^{-M}} \nu_\alpha \Bigl(
e^{(\mu_X^M-\delta)/\alpha} \vee e^{\mu_Y^M -\delta},
e^{(\mu_X^M+\delta)/\alpha} \wedge  e^{\mu_Y^M +\delta} \Bigr).
\end{align*}
\\
 Since $\mu_X^M \approx 1$ and $\mu_Y^M \approx 1/\alpha $,
  by picking $M$ large and $\delta $ small,   the right side can be
  made to be less than 1. This completes the proof.
\end{proof}

\bibliographystyle{plainnat}

\def\cprime{$'$}

\end{document}